
\input amstex
\input amssym
\documentstyle{amsppt}
\hsize = 5.4 truein
\vsize = 8.7 truein
\baselineskip = 24pt
\NoBlackBoxes
\TagsAsMath
\define\rest{\upharpoonright}
\define\forces{\Vdash}

\define\si{\sigma}

\define\la{\lambda}
\define\k{\kappa}

\define\a{\alpha}
\define\be{\beta}
\define\de{\delta}

\define\ga{\gamma}
\define\al{\aleph}

\define\sm{\setminus}

\define\th{\theta}
\define\N{\Cal N}

\define\om{\omega}

\define\orh{\overset{\rightharpoonup}\to}

\define\z{\zeta}
\define\h{\chi}
\define\m{\mu}

\redefine\v{\vert}
\define\noi{\noindent}
\define\1{\bigskip}
\define\se{\subseteq}
\define\es{\emptyset}
\define\ci{[i]_}
\define\cj{[j]_}
\define\3item{\par\indent\indent\hangindent2\parindent\textindent  }
\define\bcu{\bigcup}
\define\pr{\leq^{pr}_}
\define\apr{\leq^{apr}_}
\define\stm{\setminus}
\define\dm{dom\ }
\define\nin{\not\in}
\define\cd{card\ }
\define\Qapr{Q^{apr}_}
\define\weakpower{(2^{<\th_\k})}
\define\prm{\prime}
\define\sse{\subseteq}
\topmatter
\title Forcing Many Positive Polarized Partition Relations\\
Between a Cardinal and its Powerset
\endtitle
\author Saharon Shelah and Lee J. Stanley
\endauthor
\thanks {The research of the
first author was partially supported by the NSF and
the Basic Research Fund, Israel Academy of Science.  This is
paper number 608 in the first author's list of publications.
VERSION OF September 13, 1997}
\endthanks
\address{Hebrew University, Rutgers University}
\endaddress
\address{Lehigh University}
\endaddress

\abstract{We present a forcing for blowing up $2^\la$ and
making
\lq\lq many positive polarized partition relations"
(in a sense made precise in (c) of our main theorem)
hold in the
interval $[\la,\ 2^\la]$.  This generalizes results of
\cite{276}, Section 1, and the forcing is a \lq\lq many
cardinals" version of the forcing there.}
\endabstract

\endtopmatter
\subheading{\S 0.  INTRODUCTION}

In \cite{276}, the first author proved (with $\mu$ in
the place of our $\la$, and $\la$ in the place of our
$\eta$) the consistency of:

$$\la < \k < \eta \text{ are all regular, }2^\la = \eta,\ \eta \rightarrow
(\eta,[\k;\k])$$

The forcing can be thought of as a \lq\lq filtering through"
$\k$ of adding $\eta$ many Cohen subsets of $\la$.
Then, $\{\la,\ \k,\ \eta \}$ can be thought of as a
three element set $K$ of regular cardinals used for defining
the forcing; the elements of $K$ are taken, in the ground
model, to be sufficiently far apart.  An important technical
notion, related to the idea of \lq\lq filtering through"
is the possibility of viewing $p \leq q$ as split up,
in various ways, into \lq\lq pure" and \lq\lq apure" extensions.

It is natural to attempt to allow the set $K$ of regular
cardinals to be larger, and to simultaneously obtain
many such, and stronger, partition relations, for example, by
increasing the \lq\lq dimension"
(from 2 to $n$) and the number of blocks (from 2 to $\si$).
These will all be aspects of our treatment here, see (B), below,
and (c) of our main Theorem.

More specifically,
we start, in $V$, from

$$\text{(A)  } cf\ \la = \la = \la^{<\la} < \mu = \mu^\la = cf\ \mu,$$

\noi
and we fix

$$\text{(B)  } K \se [\la,\ \mu]\text{, a set of regular cardinals, with }
\la,\ \mu \in K.$$

In \S 1, we define a forcing $\bold Q = \bold Q_K$ which
generalizes the forcing of \cite{276}, \S 1, and we prove
its important properties, culminating in (1.13) and (1.14),
whose statements are incorporated into our main Theorem, below
(everything except item (c)).

Item (c) of the main theorem, below,
which we address in \S 2, and
which speaks of the positive polarized partition
relations, could be vacuous,
unless the elements of $K$ are sufficiently far apart, viz.
Remark 1, below.
However, nothing about the
forcing depends on the \lq\lq spacing" of elements of $K$ so no such
assumptions about $K$ figure as hypotheses until item (c).

For the positive polarized partition relations,
in $V^{\bold Q}$,
we suppose $\si,\ \k,\ \k_1,\ \h,\ \tau,\ \k_2$
are cardinals satisfying:

$$\text{(C)  } \si < \la,\ \k \leq \k_1 \leq \h = \h^\si < \tau \leq
\k_2 \text{, where } \k_1,\ \k_2 \text{ are successive members of } K.$$

If $2 \leq n < \om$, then, examining the methods of \cite{289}:

$$\text{(D)  there is } m = m(n) < \om
\text{ sufficiently large that there is a system
as in (2.1), below.}$$

Note that $m$ depends only on $n$ and not
on the cardinals of (C).  This justifies our notation.
We can now state our main theorem.
\1

\proclaim{Theorem}  If, in $V,\ \la,\ \mu,\ K$ are as in (A), (B),
above,
then there is $\bold Q = \bold Q_K= (Q,\ \leq )$ such that
the empty condition of $\bold Q \forces
2^\la \geq \mu$ and forcing with $\bold Q$
adds no sequences of length $<\ \la$.
Further, assuming that in
$V,\ 2^\th = \th^+$ for all cardinals $\th \in [\la,\ \mu]$:

\roster
\item"{(a)}"  $card\ Q = \mu$.
\item"{(b)}"  Forcing with $\bold Q$ preserves cofinalities, and therefore
cardinals.
\item"{(c)}"  Suppose that the cardinals $\si,\ \k,\ \k_1,\ \h,\ \tau,\ \k_2$
are as in (C), above, let $2 \leq n < \om$ and let $m = m(n)$ be as in (D).
If, in $V,\ \tau \rightarrow (\k)^m_\h$, then, in $V^{\bold Q},\
\ ((\tau)_\si) \rightarrow ((\k)_\si)^{((1)_n)}_\h$.
\endroster
\endproclaim
\1
\proclaim{Remarks}

\noi
\roster
\item  We remind the reader of the meaning of the second partition symbol
in the statement of the Lemma.  Let $(X_i\ :\ i < \si \}$ be a pairwise
disjoint family of sets each of cardinality (at least) $\tau$,
let $X = \bigcup\{ X_i\ :\ i < \si \}$, and
let $D$ be the set of $n-\text{element}$ subsets of $X$ which meet each
$X_i$ in at most one element.  The partition symbol then asserts that
whenever $F$ is a function from $D$ to $\h$, for $i < \si$,
there is $Y_i \se X_i$, of cardinality (at least) $\k$ such that,
letting $Y = \bigcup\{Y_i\ :\ i < \si\},\ F$
is constant on $D \cap [Y]^n$

\item  By our hypotheses
on cardinal exponentiation in $V$, it is only
the \lq\lq spacing" of the elements of $K$  that will
determine how often, in item (c) of the Theorem, we have,
in $V$, the hypothesis that $\tau \rightarrow (\k)^m_\h$.
Thus, it is only the spacing of the elements of $K$
which determines how \lq\lq many" of these positive
polarized partition relations hold in $V^\bold Q$.  Further, as
is usually the case, the assumption of GCH is just for notational
convenience and to be able to state
results in a simple compact form.  The technical lemmas of \S\S 1, 2 are
stated in a form which makes no assumptions about cardinal exponentiation,
and which indicates how the statement of the Theorem could be modified
so as not to appeal to GCH.

\item  In (B), above, we have omitted a plausible hypothesis, namely
that for $\k_2 \in K,\ sup\ K \cap \k_2 \in K$ (and so, in
particular, is regular).  This is because,
as in the proceeding Remark, the only effect of this will be to enlarge
the set of instances in which our theorem applies.  Further, by simply
enlarging the given set $K$ in the obvious ways (adding the supremum, if
it is regular, and not already in $K$, or adding the successor of the
supremum, if the supremum is singular), we can achieve the effect of this
hypothesis.

\item  In (C), above, given the role of $\tau$, we can, of course,
allow $\tau > \k_2$, but the interesting case is when
$\tau \leq \k_2$, and in fact, when $\tau < \k_2$.  Nevertheless,
the case $\tau = \k_2$ can also be handled and we indicate how to
do this at the end of (2.2).

\item  Regarding $\k$, clearly the most interesting case is
when $\k = \k_1$; unfortunately at this point, it is unclear
whether our methods, or a small variant thereof will suffice
to handle this case.  We are continuing to investigate this question
and also the question of whether we can allow $\si = \la$, at least
under the additional assumption that $\la$ is not strongly inaccessible.

\item  In order to handle all $n < \om$ simultaneously,
it is natural to use measurable cardinals and
and the obvious attempt to do so works in a straightforward.  Some
significant use of large cardinals is necessary.

\item  We treat only the extremely dispersed case, where,
in the n-tuples in the domain, each coordinate comes from a
different one of the $\si$ many blocks (the superscript $((1)_n)$).
It would be very desirable to allow pairs, or more, from the
same block.  This paper does not address this question, but
for one pair, see \cite{276}, \S 2, \cite{288}, \cite{346},
\cite{481} and \cite{585}.

\item  We began work on this paper in 1986, using essentially
the same approach as presented here; this work has been subject to
various interruptions which has made us decide to finally present it in its
present form rather than attempt to polish off various of the small questions
indicated above and to optimize the results.

\item  Our notation and terminology is intended to either
be standard, have a clear meaning, or be specifically
introduced, as needed.
\endroster
\endproclaim
\1\1
\subheading {\S 1.  THE FORCING}

We present the forcing $\bold Q$ and develop its
basic properties.  As mentioned
above, $\bold Q$ is a \lq\lq many cardinals" generalization of the forcing
of \cite{276}, \S 1.

\proclaim{(1.1)  Context and Preliminaries}
\endproclaim

Let $\la = \la^{<\la},\ \m = \m^\la,\ \la,\ \m$ both be regular.  Let $K \se
[\la,\ \m]$ be a set of regular cardinals with $\la,\ \m \in K$.  For
the remainder of this paper, $\la,\ \mu,\ K$ are fixed.

For $\k \in K$, let $E_\k$ be the equivalence relation on $\m$ defined by
$i\ E_\k\ j$ iff $i\ +\ \k = j\ +\ \k$.  For $\la \leq \k \leq \m$, define
$E_{<\k}$ as $id_\m \cup \bigcup \{E_\th:\ \th \in K \cap \k \}$.  For such
$\k$, if $\k \not\in K$, let $E_\k = E_{<\k}$.  For $i < \m,\ \la \leq \k \leq
\m$, let $[i]_\k =\ $ the $E_\k-$equivalence class of $i$, and for $A \se \m$,
let $A/E_\k = \{[i]_\k:\ i \in A \}$.  For such $i,\ A,\ [i]_\k$ is \bf
represented in $\bold{A}$ \rm iff $A \cap [i]_\k \ne \es$.
If $A \se B \se \m$, the $\ci\k$ \bf grows from
$\bold A$ to $\bold B$ \rm iff $\es \neq A \cap \ci\k \ne B \cap \ci\k$

\proclaim{(1.2)  Remarks}
\roster
\item  If $\th < \k$, both in $K$, then $E_\th$ refines $E_\k$ and, in fact,
each $E_\k$ class is the union of $\k {\text\ many\ } E_\th$ classes.
\smallskip
\item  For all $i,\ j < \m,\ i\ E_\m\ j$.  Thus, the following definition makes
sense:

$$\text{if } i < j < \m,\ \k(i,\ j) =\ \text{the least } \k \in K\ \text{such
that } i\ E_\k\ j.$$

\endroster
\endproclaim

\proclaim{(1.3) Definition and Remark}

Suppose $\k \in K$.  We define $\th_\k$ to be the least regular cardinal
which is $\geq \bigcup[(K \cap \k) \cup \{\la\}]$.  Thus, in particular,
$\th_\la = \la$, if $\th < \k$ are successive elements of $K$ then $\th_\k =
\th$, if $\bigcup (K \cap \k)$ is singular, then $\th_\k = (\bigcup (K \cap
\k))^+$, while if $\bigcup (K \cap \k)$ is inaccessible, then $\th_\k = \bigcup
(K \cap \k)$.
\endproclaim

\proclaim{(1.4)  Definition}  $q \in Q = Q_K$ iff $q:\ dom\ q\  \rightarrow
\ \{0, 1 \},\ dom\ q \se \m$ and:

\roster
\item"{(a)}"  for $i <\m,\ \k \in K,\ card (\ci\k \cap dom\ q) < \th_\k$ (note:
taking $\k = \m$, we have $card\ dom\ q < \th_\m$).
\endroster

If $p,\ q \in Q$, we set $p \leq q$ iff

\roster
\item"{(b)}"  $p \se q$,

\item"{(c)}"  For all $\k \in K,\ \{A \in \m/E_\k: A$ grows from $dom\ p$
to $dom\ q \}$ has power $< \th_\k$.
\endroster

$\bold Q = (Q,\ \leq)$.
\endproclaim

\proclaim{(1.5)  Remarks}
\roster
\item The content of (a) of (1.4) is that not too
much (less than $\theta_\kappa$) of any $E_\kappa$ class is present in
the domain of any condition.
\smallskip
\item The content of (c) is that
few (less than $\theta_\kappa$)
$E_\kappa$ classes grow from $dom\ p$ to $dom\ q$.
\item It should be emphasized that the definition of \lq\lq $A$
grows from $dom\ p$ to $dom\ q$ {\bf requires that } $A \cap dom\ p
\neq \es$.
\endroster
\endproclaim

\proclaim{(1.6)  Definition}  For $\k \in K$ and $p,\ q \in Q$, let:
$p \leq^{pr}_\k q$ iff $p \leq q$ and:

\roster
\item"{(d)}"  no $E_\k-$class represented in $dom\ p$ grows from $dom\ p$
to $dom\ q$,
\endroster

and let:  $p \leq^{apr}_\k q$ iff $p \leq q$ and:

\roster
\item"{(e)}"  $(dom\ q)/E_\k = (dom\ p)/E_\k$.
\endroster
\endproclaim

\proclaim{(1.7)  Remarks}
\roster
\item The content of (d) of (1.6) is that the only
new elements of $dom\ q$ are in new $E_\kappa$ classes.
\smallskip
\item The content
of (e) of (1.6) is that no new $E_\kappa$ classes
are represented in $dom\ q$, or,
seen in another light, that if $i \in dom\ q \setminus dom\ p$, then
$\ci\k$ grows from $dom\ p$ to $dom\ q$.
\endroster
\endproclaim

\proclaim{(1.8)  Proposition}
\roster
\item"{(a)}"  For all $\k \in K,\ \pr\k,\ \apr\k$ are partial orderings of
$Q$.

\item"{(b)}"  If $p_1,\ p_2 \in Q$ and they are compatible as functions,
then $p_1 \cup p_2 \in Q$; further, letting $q = p_1 \cup p_2$, if
(c) of (1.4) holds between $p_i$ and $q$, for $i = 1,\ 2$, then $q$
is the join, in $\bold Q$, of $p_1$ and $p_2$.

\item"{(c)}"  If $p \leq q,\ \k \in K$, then there are $r,\ s \in Q$ such
that:
\itemitem  {(1)}  $p \leq^{pr}_\k\ r \leq^{apr}_\k q$,
\itemitem  {(2)}  $p \leq^{apr}_\k s \leq^{pr}_\k q$ and
\itemitem  {(3)}  $q = r \cup s$.

\item"{(d)}"  $\leq = \leq^{apr}_\m$ (except that if $\es \neq q \in Q$,
then $\es \leq q$, but $\es \not\leq^{apr}_\k q$ for any $\k \in K$).

\item"{(e)}"  If $\k_0 \leq \k_1 \leq \k_2$, all $\in K$, then:

\noindent $\pr{\k_1} \se \pr{\k_0},\ \apr{\k_1} \se \apr{\k_2}$.

\item"{(f)}"  If ($\k \in K\ \&\ s\ \apr\k\ t\ \&\ s\ \pr\k\ v$), then
$t \cup v \in Q$ and:

\noindent $s \leq (t \cup v),\ t\ \pr\k\ (t \cup v),\ v\ \apr\k\ (t \cup v)$.

\item"{(g)}"  If $\k \in K,\ p\ \leq^*_\k\ q_i\ (i = 1,\ 2)$, where $* \in
\{pr,\ apr \}$ and $q_1,\ q_2$ are compatible in $\bold Q$, then
$p\ \leq^*_\k\ (q_1 \cup q_2)$.

\item"{(h)}"  If $p \leq^{apr}_\k q_1,\ q_2$ and if

$$\text{(*) if } (i \in dom\ q_1 \setminus
dom\ p\ \&\ j \in dom\ q_2 \setminus dom\ p) \text{ then }(\ci{<\k}
\neq \cj{<\k} \text{ or } \ci{<\k} \cap dom\ q_1 = \cj{<\k} \cap dom\ q_2))$$

then also $q_k \leq^*_\k q_1 \cup q_2,\ k = 1,\ 2$.

\item"{(i)}"  If $p \leq^{apr}_\k\ q_i \leq r$ for $i = 1,\ 2$, then,
for such $i,\ q_i \leq^{apr}_\k q_1 \cup q_2$.
\endroster
\endproclaim

\noindent
\demo{Proof}  (a) and (b) are clear.
For (c), let $r = q \v x$, where
$\xi \in x$ iff $\xi \in dom\ q$ and $(\xi \in dom\ p$ or $[\xi]_\k \cap
dom\ p = \es)$.  Also, let $s = q \v y$, where $\xi \in y$ iff $\xi \in dom\ q$
and $[\xi]_\k \cap dom\ p \neq \es$.  Clearly $p\ \pr\k\ r,\ p\ \apr\k\ s$;
clearly $q = r \cup s$.
We verify that $r\ \apr\k\ q$ and $s\ \pr\k\ q$.  For the first, suppose that
$\xi \in dom\ q \setminus x$.  Then, $\xi \not\in dom\ p$ and
$[\xi]_\k \cap dom\ p \neq \es$.  Then certainly $[\xi]_\k \cap x \neq \es$,
i.e. $\xi \in \bcu dom\ r/E_\k$.  For the second, suppose
$\xi \in dom\ q \setminus y$, but $[\xi]_\k \cap y \neq \es$.  Then,
$\xi \in x \setminus dom\ p$, so $[\xi]_\k \cap dom\ p = \es$.  If
$\z \in [\xi]_\k \cap y$, then $[\z]_\k \cap dom\ p \neq \es$, but
$[\z]_\k = [\xi]_\k$, contradiction.

For (d), recall that $\m = \ci\m$, for all $i < \m$.  For (e),
if $p\ \pr{\k_1}\ q$ and $x$ is an $E_{\k_0}-$class represented in $p$,
let $x^*$ be the $E_{\k_1}-$class such that $x \se x^*$.  Then, $x^*$ is
represented in p and since $x^*$ does not grow from $dom\ p$ to $dom\ q$,
neither can $x$.  Similarly, if $p\ \apr{\k_1} r$ and $\xi \in dom\ r$,
there is $\z\ E_{\k_1}\ \xi$ such that $\z \in dom\ p$.  But then,
$\z\ E_{\k_2}\ \xi$, so $\xi \in [\z]_{\k_2}$.

For (f), we first show that $(\dm t \stm \dm s) \cap (\dm v \stm \dm s) =
\es$; then, by (b), $t \cup v \in Q$.
It will then be clear that $s \leq (t \cup v)$.  So,
if $\xi \in \dm v \stm \dm s$, then $[\xi]_\k \cap \dm s = \es$, so
$\xi \nin \bcu \dm s/E_\k$ and $\dm t \se \bcu \dm s/E_\k$.

Next, we show that $t,\ v \leq t \cup v$; by (b), it will suffice to
show that (c) of (1.4) holds between $t$ and $t \cup v$ and between
$v$ and $t \cup v$.  We prove the former first.
So, suppose that $\tau \in K$ and
first suppose that $\es \neq dom\ t \cap \ci\tau$ and
$j \in (dom\ v \cap \ci\tau) \setminus dom\ t$.  Then, certainly $j \not\in
dom\ s$, so since $s \leq^{pr}_\k v$, we clearly must have that $\tau > \k$.
Now, let $l \in \ci\tau \cap dom\ t$.  Since $s \leq^{apr}_\k t$,
there is $a \in dom\ s$ such that $a\ E_\k\ l$.  But then, since
$\tau > \k$ (actually, $\geq$ would suffice here), $a\ E_\tau\ i$, so
$\es \neq dom\ s \cap \ci\tau \neq dom\ s \cap \ci\tau$.
And, since $s \leq v$, there are fewer than $\theta_\tau$ many such $\ci\tau$,
and we have proved that (c) of (1.4) holds between $t$ and $t \cup v$.

To show that (c) of (1.4) holds between $v$ and $t \cup v$, let $\tau$
be as above, and, this time, suppose that $\es \neq  dom\ v \cap \ci\tau$
and that $j \in (dom\ t \cap \ci\tau) \setminus dom\ v$.  Then, certainly
$j \not\in dom\ s$, and so, since $s \leq^{apr}_\k,\
\cj\k \cap dom\ s \neq \es$.  Thus, $\cj\k$ grows from
$dom\ s$ to $dom\ t$, and, since $s \leq t$, there are at most $\theta_\k$
many such $\cj\k$.  We consider separately the cases $\tau \geq \k$
and $\tau < \k$.  In the first case,
$\theta_\k \leq \theta_\tau$ and we have found
one of at most $\theta_\k$ many $\cj\k$ inside every $\ci\tau$
which grows from $dom\ v$ to
$dom\ t \cup v$, so clearly there are at most $\theta_\tau$ many
such $\ci\tau$, as required.
Thus, without loss of generality, we may assume that
$\tau < \k$.  In this case, we shall argue that
$\es \neq dom\ s \cap \ci\tau$.  Clearly this will suffice since
then $\ci\tau$ grows from $dom\ s$ to $dom\ t$, and again, since
$s \leq t$, there are at most $\theta_\tau$ such $\ci\tau$, as
required.

So, suppose, towards a contradiction, that $\es = dom\ s \cap \ci\tau$.
Let $\xi \in \ci\tau \cap dom\ v$, so $[\xi]_\k = \ci\k$.  But
$j \in \ci\tau \cap dom\ t$, so $j \in \ci\k \cap dom\ t$.  Since
$s \leq^{apr}_\k t$, this means that $\es \neq \ci\k \cap dom\ s$.
But then $\xi \in \ci\k \cap (dom\ v \setminus dom\ s)$.  This,
however, is impossible, since $s \leq^{pr}_\k v$, which completes
the proof.

We proceed, now to show that $t \leq^{pr}_\k t \cup v$ and that
$v \leq^{apr}_\k t \cup v$.  For the former, suppose that
$\xi \in \dm v \stm \dm t$.  Then, $\xi \in \dm v \stm \dm s$, so
$[\xi]_\k \cap \dm s = \es$.  We claim that $[\xi]_\k \cap \dm t = \es$.
If not, and $\z \in [\xi]_\k \cap \dm t$, then $[\z]_\k \cap \dm s \neq \es$,
but, once again, $[\z]_\k = [\xi]_\k$, contradiction.  Thus,
$t\ \pr\k\ (t \cup v)$.

To see that $v \leq^{apr}_\k t \cup v$,
suppose that $\xi \in \dm t \stm \dm v$.  We
need to show that $[\xi]_\k \cap \dm v \neq \es$.
This, however, is clear, because, since
$\xi \in \dm t,\ [\xi]_\k \cap \dm s \neq \es$, so certainly
$[\xi]_\k \cap \dm v \neq \es$, and we have finished proving (f).

For (g), first note that if $q_i \leq r$ for $i = 1,\ 2$, then,
letting $s = q_1 \cup q_2$,
for such $i,\ q_i \leq s \leq r$.  This is clear,
because if $\tau \in K$ and $\cj\tau$ grows from $dom\ q_i$
to $dom\ s$, then certainly $\cj\tau$ grows from $dom\ q_i$
to $dom\ r$, and there are at most $\theta_\tau$ such $\cj\tau$,
since $q_i \leq s$.  Further, if $\cj\tau$ grows from
$dom\ s$ to $dom\ r$, then either $\cj\tau$ grows from
$dom\ q_1$ to $dom\ r$ or $\cj\tau$ grows from $dom\ q_2$
to $dom\ r$, and again, since $q_1,\ q_2 \leq r$, there
are at most $\theta_\tau$ such $\cj\tau$ for each case.

Now suppose that $*$ is $apr$.  Thus, if $\xi \in \dm s$,
then, for an $i \in \{1,\ 2\},\ \xi \in \dm q_i$, so
$[\xi]_\k \cap \dm p \neq \es$.  It is then clear that $p \leq^*_\k s$,
as required.

If $*$ is $pr$ and
$\xi \in \dm s \stm \dm p$, then, letting
$i \in \{1,\ 2\}$ be such that $\xi \in \dm q_i \stm \dm p$, then, since
$p\ \pr\k\ q_i$, clearly $[\xi]_\k \cap \dm p = \es$, as required.

We prove (i), before proving (h).  As in (g), let
$s = q_1 \cup q_2$.
For $i = 1,\ 2$, we must show
that $q_i \leq^{apr}_\k s$.
We already know, from the proof of (g), that for such $i,\
q_i \leq s$.  So, let
$j = 1 + (2 - i)$,
and suppose that $\a \in dom\ s \setminus
dom\ q_i$.  We need to show that $\es \neq [\a]_\k \cap dom\ q_i$.
But $\a \in dom\ q_j \setminus dom\ q_i$, so $\a \in dom\ q_j
\setminus dom\ p$, so $\es \neq [\a]_\k \cap dom\ p$, and the
conclusion is clear.

We conclude by proving (h).  For this,
let $s = q_1 \cup q_2$.  If we prove that $q_1$ and $q_2$
are compatible in $\bold Q$, then, by (g) and (i), we are finished.
In fact, we will show directly that $q_1,\ q_2 \leq s$.
By symmetry, it will suffice
to prove that $q_1 \leq s$, and clearly, only (c) of (1.4) is at issue.
So, let $\tau \in K$.  First note that, without loss of generality,
we may assume that $\tau < \k$.
This is because, since $dom\ q_1 \stm dom\ p$ and
$dom\ q_2 \stm dom\ p$ both have cardinality less than $\theta_\k$,
therefore so do $dom\ q_1 \stm dom\ q_2$ and $dom\ q_2 \stm dom\ q_1$.
Then, if $\tau \geq \k$,
in particular, fewer than $\theta_\k$ many $E_\tau$ classes grow
from $dom\ q_1$ to $dom\ s$.

So, suppose $\tau > \k$.
By hypothesis, if $\ci\tau$ grows from $dom\ q_1$ to
$dom\ s$, then $\ci\tau \cap dom\ q_1 = \ci\tau \cap dom\ p \neq \es$,
and so $\ci\tau$ grows from $dom\ p$ to $dom\ q_2$.
However, since $p \leq q_2$, there are fewer than
$\theta_\tau$ such $\ci\tau$.  This concludes the proof of (h) and
of the Proposition.

\enddemo

\proclaim{(1.9)  Proposition}
\roster
\item"{(a)}"  For all $\k \in K,\ (Q,\ \pr\k)$ is
$\k\text{-complete.}$

\item"{(b)}"  $\bold Q$ is $\la\text{-complete.}$
\endroster
\endproclaim

\noindent
\demo{Proof}  For (a), let $\th < \k$ be regular and let $(q_i: i < \th)$ be a
$\pr\k-\text{increasing}$ sequence.  We claim that $q =_{def}
\bcu \{q_i: i < \th \} \in Q$ and that for all $i < \th,\ q_i\ \pr\k\ q$.
We first verify (a) of (1.4).  For this, let $i < \m,\ \nu \in K$.  Clearly
if $\k < \nu$, then $\cd (\ci\nu \cap \dm q) < \th_\nu$, since
$\th < \k \leq \th_\nu,\ \th_\nu$ is regular, $\ci\nu \cap \dm q =
\bcu \{\ci\nu \cap \dm q_j: j < \th \}$ and for all $j < \th,
\ \cd (\ci\nu \cap \dm q_j) < \th_\nu$.  So, suppose $\nu \leq \k$.  If
$\ci\nu \cap \dm q = \es$, there is nothing to prove, so suppose that
$j < \th$ and $\xi \in \ci\nu \cap \dm q_j$.  Thus,
$\ci\nu \cap \dm q_j \neq \es$, and therefore $\ci\k \cap \dm q_j \neq \es$.
But then, for $j < \ell < \th$, since $q_j\ \pr\k\ q_\ell,
\ \ci\k \cap \dm q_\ell = \ci\k \cap \dm q_j$, so, in fact,
$\ci\k \cap \dm q = \ci\k \cap \dm q_j$ and therefore
$\ci\nu \cap \dm q = \ci\nu \cap \dm q_j$, which yields (a).

Clearly $q_i \se q$.  Further, if we verify that \bf no \rm
$E_\k-\text{class}$ represented in $\dm q_i$ grows from
$\dm q_i$ to $\dm q$, then
of course we have verified (c) of (1.4) and so $q_i \pr\k q$.  So, let $A$ be
an $E_\k-\text{class}$ represented in $\dm q_i$.  If $\xi \in A \cap \dm q$,
then, for some $i < j < \th,\ \xi \in A \cap q_j$, so, since
$q_i\ \pr\k\ q_j,\ \xi \in \dm q_i$, and we are finished.

For (b), let $\th < \la = \theta_\la$, and let
$(q_i:\ i < \th)$ be increasing for $\leq$.  Let $q =
\bigcup \{q_i:\ i < \th\}$.  We must verify (a) of (1.4) for
$q$, and that for all $i < \th$, (c) of (1.4) holds between $q_i$
and $q$.  So, let $\k \in K$, and let $a$ be an
$E_\k\text{-class}$.  Clearly $a \cap dom\ q =
\bigcup \{a \cap dom\ q_i:\ i < \th\}$ and since each $q_i \in Q$,
each $a \cap dom\ q_i$ has power less than $\theta_\k$, and therefore
their union also has power less than $\theta_\k$, since $\theta_\$
\geq \la$ and $\theta_\k$ is regular.  Similarly, if $i < \th$,
then $\{a:\ a\text{\ is an\ } E_\k\text{-class which grows from\ }
dom\ q_i\text{\ to\ } dom\ q\} =$

$\noindent \bigcup \{X_j:\ i < j < \th \}$, where,
for $i < j < \th,\ X_j\ :=\ $

$\noindent\{a:\ a\text{ is an } E_\k\text{-class which grows from }
dom\ q_i\text{ to } dom\ q_j\}$.  Again, since, for $i < j < \th,\
q_i \leq q_j$, by (c) of (1.4), each $X_j$ has power less than
$\theta_\k$, and therefore, as before, the same is true of their
union.  This completes the proof of the Proposition.
\enddemo

\proclaim{(1.10)  Proposition}  If $\k \in K,\ p \in Q$, then
$\bold {\Qapr{\k,\ p}}$ has the $(2^{<\th_\k})^+-\text{c.c.}$, where
$\bold {\Qapr{\k,\ p}} = (\Qapr{\k,\ p},\ \apr\k)$, and
$\Qapr{\k,\ p} = \{q: p\ \apr\k q \}$.
\endproclaim

\noindent
\demo{Proof}  We should note, here, immediately, that
in virtue of (1.8), (i), for $q_1,\ q_2 \geq p$,
compatibility in $\bold {\Qapr{\k,\ p}}$
is the same as compatibility in $\bold Q$, so it
is the latter that we shall establish, when
our statement calls for the former.

Suppose, now, that $q_i \in \Qapr{\k,\ p}$, for
$i < \weakpower^+$.  We show there is $I \se \weakpower^+$ with
$\cd Y = \weakpower^+$,
such that for $i,\ j \in I,\ q_i$ and $q_j$ are compatible
in $\bold Q$.  In virtue of the preceding paragraph,
clearly this suffices.

For $i < \weakpower^+$, let
$d_i = \dm q_i \stm \dm p$.  We first show
that $\cd d_i < \th_\k$.  Note that by (e) of (1.6),
if $\a \in d_i$, then $[\a]_\k$ grows from $dom\ p$ to
$dom\ q$, and so $d_i/E_\k \se
\{ A \in \mu(\k):\ A \text{ grows from } dom\ p \text{ to } dom\ q\}$.
By (1.4), (c), this last set has power $< \theta_\k$.  Finally,
by (1.4), (a), for all $A \in d_i/E_\k,\ \cd (A \cap dom\ q_i) < \theta_\k$.
Then, since $\theta_\k$ is regular, the conclusion
that $card\ d_i < \th_\k$ is clear.

Consider now $Y_i := d_i/E_{<\k}$.  Since each $Y_i$ has power $< \theta_\k$,
by the $\Delta\text{-system}$ Lemma,
there is $Y \se \mu /E_{<\k}$ and $I \se \weakpower^+$,
with $\cd Y = \weakpower^+ < \theta_\k$ such that for $i,\ j \in I,\ Y_i \cap
Y_j = Y$.  Now, $\cd Y < \theta_\k$, so there are at
most $\theta_\k^{<\theta_\k} = \weakpower$ many possible
$d_i \cap \bigcup Y$, and therefore, without loss of generality, there is
$d$ such that for $i \in I,\ d_i \cap \bigcup Y = d$.  Note
that for this $d$, there are at most $2^{\cd d} \leq \weakpower$
many possible $f:d \rightarrow 2$, so, again, without loss of
generality, there is $f$ such that for $i \in I,\ q_i\rest
(d_i \cap \bigcup Y) =
f$.

By (1.8), (2), for $i,\ j \in I,\ q_i \cup q_j \in Q$, and in order
to conclude, as usual,
it will suffice to show that, letting $s = q_i \cup q_j$,
clause (c) of (1.4)
holds between $q_i$ and $s$ (by symmetry, the same will be true,
of course, replacing $q_i$ by $q_j$).  So, let $\tau \in K$, and
suppose that $a \in \mu/E_\tau$ and $a$ grows from $dom\ q_i$ to
$dom\ s$.  As in the proof of (1.8), (h), without loss of generality,
$\tau < \k$.  We argue that our hypotheses imply that
$\es \neq a \cap dom\ p$.  Of course, this will suffice, since
$p \leq q_j$, and so there are fewer than
$\theta_\tau$ such $a$.

So, suppose, towards a contradiction, that $\es =
a \cap dom\ p$.
But, in this case, $a \in Y_j$.  Since
$a$ grows from $dom\ q_i$ to $dom\ s,\ \es \neq a \cap dom\ q_i$,
so also $a \in Y_i$, so $a \in Y$.  But now we have a contradiction,
since by properties of $I,\ d_i \cap \bigcup Y = d =
d_j \cap \bigcup Y$, and so $a$ cannot grow from $dom\ q_i$
to $dom\ s$, after all, contrary to our hypothesis.  This completes
the proof of the Proposition.
\enddemo
\medskip

We need a slightly more refined version of this.
\medskip

\noindent
\proclaim{(1.11)  Proposition}  Suppose $\k \in K,\ \weakpower^+ \leq \k,
\ (s_i: i < i^*)$ is a $\pr\k-\text{increasing}$ sequence from $Q$, and
suppose that for $i < i^*,\ s_i\ \apr\k\ t_i$, and that for $j < i < i^*,
\ t_j,\ t_i$ are incompatible in $\bold Q$.  Then, $i^* < \weakpower^+$.
\endproclaim

\demo{Proof}  If $i^* < \k$, we can take $s = \bcu \{s_i: i < i^*\}$.  Noting
that for $j < i^*,\ s\ \apr\k\ (s \cup t_j)$, we can then apply (1.10).  Even
if $\k \leq i^*$, we can \bf essentially \rm argue in this fashion, by
redoing the proof of (1.10).  So, let $i^* = \weakpower^+ \leq \k$.  Let $d_i =
\dm t_i \stm \bcu \{\dm s_i: i < i^* \}$.  We obtain a contradiction.
Then, $d_i \se \dm t_i \stm \dm s_i$, and,
arguing as in (1.10),
$\cd d_i \se \cd (\dm t_i \stm \dm s_i) < \th_\k$.

As in (1.10), for $i < i^*$, let $Y_i = d_i(<\k)$.
Once again, we can find $I \se i^*,\ Y,\ d$, and $f$ such that
$\cd I = i^*$ and for $i,\ j \in I,\ Y_i \cap Y_j = Y,\
d_i \cap \bigcup Y = d$ and $t_i\rest d = f$.  The conclusion is
then as in (1.10) that for $i,\ j \in I,\ t_i$ and $t_j$
are compatible in $\bold Q$ and therefore in $\bold {\Qapr{\k,\ p}}$.
This contradiction completes the proof of the Proposition.
\enddemo

\proclaim{(1.12)  Lemma}  If $\k \in K,\
2^{<\th_\k} < \k,\ p \in Q$ and
$p\ \v\vdash_{\bold Q}\lq\lq \dot\a$ is an ordinal", THEN, there are
$q$ and $(r_i: i < i^*)$, all from $Q$, such that:

\noindent
\roster
\item"{(a)}"  $i^* < \weakpower^+$,
\item"{(b)}"  $p\ \pr\k\ q$,
\item"{(c)}"  $q\ \apr\k\ r_i$, for all $i < i^*$,
\item"{(d)}"  for some $\a_i,\ r_i\ \v\vdash_{\bold Q}\lq\lq \dot\a = \a_i$",
\item"{(e)}"  $\{r_i: i < i^* \}$ is predense above $q$.
\endroster
\endproclaim

\demo{Proof}  We shall obtain $q$ as $q_{i^*} = \bcu \{q_i: i < i^* \}$, where
$(q_i: i < i^*)$ is $\pr\k-\text{increasing}$, with $q_0 = p$.  We work by
recursion on $i$.  Having obtained $(q_j: j \leq i)$ and $(r^\prm_j: j < i)$
such that $(q_j: j \leq i)$ is $\pr\k-\text{increasing}$, the
$(r^\prm_j: j < i)$ are pairwise incompatible in $\bold Q$, $q_j\ \apr\k
\ r^\prm_j$ and there is $\a_j$ such that
$r^\prm_j \v\vdash_{\bold Q}\lq\lq \dot\a = \a_j$", note that we have the
following properties:

\noindent
\roster
\item"{(1)}"  for all $j < i,\ q_i\ \apr\k\ (q_i \cup r^\prm_j)$ (this is by
(f) of (1.8) with $s = q_i,\ t = r^\prm_j,\ v = q_i$),

\item"{(2)}"  so, letting $r^{\prm\prm}_j = q_i \cup r^\prm_j,
\ \{r^{\prm\prm}_j: j < i \} \se \Qapr{\k,\ q_i}$.

\endroster

If $\{r^{\prm\prm}_j: j < i \}$ is predense in $\bold {\Qapr{\k,\ q_i}}$,
then we take
$i^* = i,\ q = q_i,\ r_j = r^{\prm\prm}_j$, for $j < i$, and we stop.
Otherwise, there is $q^\prm \in \Qapr{\k,\ q_i}$ such that $q^\prm$ is
imcompatible with each $r^{\prm\prm}_j$.  Note that, in this case, we must have
that $q^\prm$ is incompatible in $\bold Q$ with each $r^\prm_j$, by (g) of
(1.8).  In this case, we shall have $i < i^*$, and we continue, so fix such
$q^\prm$ and let $q^\prm \leq r^\prm$ be such that for some $\a,
\ r^\prm \v\vdash_{\bold Q}\lq\lq\dot\a = \a$".  Applying (c) of (1.8), we get
$q_i\ \pr\k\ q^*\ \apr\k\ r^\prm$.  We let
$q_{i + 1} = q^*,\ r^\prm_i = r^\prm$.  By (g) of (1.8), the $r^\prm_j\ (j
\leq i)$ are pairwise incompatible in $\bold Q$.

If $i$ is a limit ordinal, $i < \k$ and the $(q_j: j < i),\ (r^\prm_j: j < i)$
are definied satisfying the induction hypotheses, we let $q_i = \bcu \{q_j: j
< i \}$ (so, by (1.9), $q_i \in Q$ and is the $\pr\k-\text{lub}$ of the
$q_j$).  We must now see that the process terminates at some $i^* <
\weakpower^+$.  If not, and if $\weakpower^+ < \k$, let $q = \bcu \{q_j: j <
\weakpower^+ \}$, and (using the above observations), for $j < \weakpower^+$,
let $r_j = r^\prm_j \cup q$.  Then, the $r_j$ are a pairwise incompatible
family in $\bold {\Qapr{\k,\ q}}$, contradicting (1.10).  If $\weakpower^+
\leq \k$, we proceed as in (1.11) to see that we must have $i^* <
\weakpower^+$, contradiction.  This means, in particular, that $i^* < \k$
and then we conclude by defining $q$ and the $r_j$ as in the case where
$\weakpower^+ < \k$, but everywhere replacing $\weakpower^+$ by $i^*$.
This completes the proof of the Lemma.
\enddemo

\proclaim{(1.13)  Proposition}  The empty condition of $\bold Q$ forces
$2^\la \geq \mu$.
\endproclaim

\demo{Proof}  For $i < \mu$, let $\bold {r}_i$ be the following
$\bold{Q}\text{-name:  } \{((\ga,\ k),\ p):
\ga < \la,\ k < 2,\ p \in Q\ \&\
p(\la i + \ga) = k \}$.  Since $\theta_\la = \la$, and for
$i < \mu,\ [\la i,\ \la i + \la) = [\la i]_\la$, we clearly have
that for $p \in \bold Q$, if $i_0 < i_1 < \mu, \cd A_j < \la$, for $j = 0,\ 1$,
where, for such $j,\ A_j = \{\ga < \la:\ \la i_j + \ga \in dom\ p \}$.
So, for such $p,\ i_0,\ i_1$, choosing $\ga \in \la \stm (A_0 \cup A_1)$,
and letting $q = p \cup \{(\la i_j + \ga,\ j):\ j < 2 \}$, we have
$p \leq q$ and $q \forces
\bold {r}_{i_0} \neq \bold {r}_{i_1}$, and the conclusion is
then clear.  This completes the proof of the Proposition.
\enddemo

\proclaim{(1.14)  Proposition} (Assuming that for
cardinals $\theta$, with $\la \leq \theta \leq \mu,\
2^\theta = \theta^+)$:

\roster
\item"{(a)}"  $card\ Q = \mu$.
\item"{(b)}"  Forcing with $\bold Q$ adds no sequences of length $< \la$.
\item"{(c)}"  Forcing with $\bold Q$ preserves cofinalities, and therefore
cardinals.
\endroster
\endproclaim

\demo{Proof}  (a) is clear, and
(b) follows easily, from (1.9), (b).  For (c), assume, towards a
contradiction, that $\tau < \si$, where both are regular, but that for
some $q \in Q,\ q \forces cf\ \si = \tau$.  By (b), we may assume
that $\tau \geq \la$.  Note that by (1.8), (d) and (1.10), with $p = \es,\
\bold Q$ has the $(2^{<\theta_\mu})^+\text{-c.c}$.  Further, under
our additional hypotheses on cardinal exponentiation, $(2^{<\theta_\mu})
\leq \mu$, so, clearly we cannot have $\si > \mu$. But then there must be
$\k \in K$ such that $\theta_\k \leq \tau < \k$.  Suppose, now, that
the $\bold Q\text{-name } \bold f$ is such that $q \forces
\bold f$ is monotone-increasing, maps $\tau$ to $\si$ and has range
cofinal in $\si$.  By (1.9), (a) and (1.12), applying (1.12) repeatedly
to each of the names $\bold f (\a)$, for $\a < \tau$, we reach a
contradiction, also using that $\weakpower \leq \tau$.
This completes the proof of the Proposition.
\enddemo

\proclaim{Remark}  Thus, in order to complete the
proof of the main Theorem of the Introduction,
it remains \lq\lq only" to prove item (c).
\endproclaim
\bigskip\bigskip

\subheading {\S 2.  THE PARTITION RELATIONS}
\medskip

In this section, we address item (c) of the main theorem of the
Introduction, in the case where $\si < \la,\ \k < \k_1$, and, at first, under
the additional simplifying assumption that $\tau < \k_2$.
For convenience, we recall the context,
and restate (c) as a Lemma, with these
additional assumptions.  The remainder of the section
will be devoted to proving this Lemma.  After the
proof is given, we will briefly indicate the small
changes necessary to accomodate the case $\tau = \k_2$.
So, let $\k_1,\ \k_2$ be successive members of $K$, let $\k < \k_1 \leq
\h = \h^\si < \tau < \k_2$, let $\si < \la$.
Assume that $2^{< \k_1} \leq \tau$ (in the context of (c) of
the main Theorem, this will follow from the Theorem's
hypotheses on cardinal exponentiation).  Recall that
for all $2 \leq n < \om$, by examination of the
methods of \cite{289}, there is
sufficiently large $m(n) < \omega$
such that, assuming that, in $V,\ \tau \rightarrow (\k_1)^{m(n)}_\h$,
then, also in $V$, there is a system as in (2.1) below.

\proclaim{Lemma}  For $2 \leq n < \om$, if,
in $V,\ \tau \rightarrow (\k)^{m(n)}_\h$,
then, in $V^{\bold Q},\
((\tau)_\si) \rightarrow ((\k)_\si)^{((1)_n)}_\h$.
\endproclaim

\proclaim{Remark}  Of course, the Lemma immediately gives
(c) of the main Theorem of the Introduction, and thus
completes the proof of the Theorem
in the cases indicated above.

We prove the Lemma
in the remainder of the section.
\endproclaim

Let
$(A_i: i < \si)$ be a sequence of sets of ordinals, each of order-type $\tau$,
such that for $i < j < \si,\ A_i < A_j$.
Let $A\ :=\ \bcu\{A_i: i < \si \}$.  Let
$D\ :=\ \{a \in [A]^n:\ card(a \cap A_i) \leq 1\text{, for all } i < \si \}$.
We often view the elements of $D$ as n-tuples, enumerated in their
increasing order.
Let $\bold c$ be a
$\bold Q\text{-name}$ for a function from $D$ to $\h$.

Let $p \in Q$.
Using the methods of \S 1, we
can find a $\pr {\k_2}\text{-increasing}$ sequence from $Q,
\ \orh{p} = (p_j: j < \eta)$,
with the following properties:

\roster
\item"{(1)}"  $\eta \leq \tau$, and $p_0 = p$,

\item"{(2)}"  for each $\orh{\a} = (\a_1,\ \cdots,\ \a_n) \in D$, there is
$j = j_{\orh{\a}}$ such that in $\Qapr{\k_2,\ p_{j+1}}$, there is a predense
set, $I_{\orh{\a}}$ of conditions deciding $\bold c (\orh{\a})$.
\endroster

\subheading {(2.1)  The system of \cite{289}}
\medskip
Now, let $\nu^*$ be a sufficiently large regular cardinal.
Fix $<^*$, a well-ordering of $H_{\nu^*}$.
For sequences $(X_t: t \in I)$, let
$u \in J(X_t: t \in I)$ iff $u \se \bcu\{X_t: t \in I \},\ \cd u \leq n$ and
for all $t \in I,\ \cd (X_t \cap u) \leq 1$.  If $u,\ v \in J(X_t: t \in I)$,
we set $u \sim v$ iff for all $t \in I,\ \cd (X_t \cap u) = \cd (X_t \cap
v)$.  By \cite{289} (and our choice of $m(n)$), we have the following.
\1\1
\proclaim{Proposition}  There are $B_i \in
[A_i]^\k\ (i < \si),\ \N_u\ (u \in J(B_i: i < \si)$, and
$h_{u,\ v}\ (u,\ v \in J(B_i: i < \si),\ u \sim  v)$ satisfying:

\roster
\item"{(3)}"  each $\N_u \prec
(H_{\nu^*},\ \in,\ <^*,\ \bold c ,\ \orh{p},\ \h,\ \k_1,\ \k_2,\ K,\ \bold Q ,
\ (I_{\orh{\a}}: \a \in D),\ \cdots )$,

\item"{(4)}"  letting $N_u = \v\N_u\v,\ N_u \cap (\bcu \{B_i: i < \si \}) =
u,\ \h \se N_u,\ \cd N_u = \h,\ N_u^{<\h} \se N_u$,

\item"{(5)}"  $N_u \cap N_v \se N_{u \cap v}$ (large cardinals
are required for $=$ in place of $\se$),

\item"{(6)}"  $(h_{u,\ v}:\ u,\ v \in J(B_i: i < \si),\ u \sim v)$
is a commutative system of isomorphisms, $h_{u,\ v}:\N_u \rightarrow
\N_v$,

\item"{(7)}"  If $u_k \sim v_k$, for $k = 1,\ 2$, then $h_{u_1,\ v_1}$
and $h_{u_2,\ v_2}$ are compatible functions, when both are defined.
\endroster
\endproclaim

\subheading {(2.2)  Completing the Proof}

In this subsection, we
complete the proof of the Lemma.
Note that our hypothesis that $\tau < \k_2$
guarantees that $\eta < \k_2$.
This is the only use we make of the hypothesis that $\tau < \k_2$.

So, let $p^* = \bcu \{p_j: j < \eta \}$.  Then, since
here, we have that $\eta < \k_2,\
p^* \in Q$, and is the $\leq^{pr}_{\k_2}$ least upper
bound of the $p_j$, by (1.9), (a).
Also, let $\ga_i = min\ B_i\ (i < \si)$, and let
$J\ :=\ J(\{\ga_i\}:\ i < \si),\ \tilde J\ :=\
J(B_i:\ i < \si)$.

\proclaim{Claim 1}  If $q \in Q^{apr}_{\k_2,\ p^*},\ u \in J$,
and $q \in N_u$, then $(dom\ q) \sm (dom\ p^*) \se N_u$.
\endproclaim
\demo{Proof of Claim 1}  $(dom\ q) \sm (dom\ p^*) \in N_u$ and it
has power $< \th_{\k_2} = \k_1 \leq \h \se N_u$, so the conclusion
is clear.
\enddemo

\proclaim{Claim 2}  There is $r \in Q,\ p^* \apr{\k_2} r$ such that:

\roster
\item  $\dm r \stm \dm p^* \se \bcu \{N_u:\ u \in J\}$
\item  for all $u \in J,\ p^* \cup (r\v N_u) \in N_u$; if further,
$\cd u = n$, then $p^* \cup r\rest N_u$ decides the value of $\bold c (u)$.
\endroster
\endproclaim
\demo{Proof of Claim 2}  Note that for the first
part of (2), it suffices to have $r\rest N_u \in N_u$,
since $p^* \in N_\es$.
Note, also, that $J$ has power
$\si$, and so we enumerate $J$ as $(u_j:\ j < \si)$.  We
shall define by recursion on $j \leq \si$
a sequence $(r_j:\ j \leq \si)$ with $r_0\ :=\ p^*$,
and all $r_j \in Q^{apr}_{\k_2,\ p^*}$.
We shall have $r\ :=\ r_\si$.
The following induction hypotheses will be in vigor,
for $j \leq \si$.  The parallel with items (1) and (2)
in the statement of the Claim should be clear.
\1
\roster
\item"{(a)}"  if $k + 1 \leq j$ then $\dm r_{k+1} \stm \dm r_k \se
N_{u_k}$,
\item"{(b)}"  for all $u \in J$ and all $k \leq j,\ r_k\rest N_u \in N_u$;
if, further, $k + 1 \leq j$ and $card\ u_k = n$,
then $p^* \cup (r_{k+1}\rest N_{u_k})$ decides the value of $\bold c (u_k)$
\item"{(c)}"  $(r_k:\ k \leq j)$ is $\apr{\k_2}$ increasing.
\endroster
\1
Clearly (a) - (c) hold for $j = 0$ with $r_0 = p^*$.
At limit ordinals, $\de \leq \si$, we shall take
$r_\de\:=\ \bigcup\{r_j:\ j < \de\}$.  If $\de < \si$,
then $\de < \la \leq \k_1 = \th_{\k_2}$.
Thus, if $\de < \si$, by (1.9), (b), $r_\de \in Q$
and is the $\leq$ least upper bound of the $r_j$.
Then, clearly also it is the $\apr{\k_2}$
least upper bound of the $r_j$.

If $\de = \si$, then, since we
are assuming $\si < \la$, we also have $\de < \la$,
and so, the same arguments yield the same conclusions,
in this case as well.

Clearly this preserves (a), (c) and the second part of (b).
We argue that it also preserves the first part of (b).
So, let $u \in J$.  We must see that $r_\de\rest N_u \in N_u$.
But $r_\de\rest N_u = \bigcup \{r_k\rest N_u:\ k < \de \}$, and
for all $k < \de$, by (the first part of) (b) for (k),
$r_k\rest N_u \in N_u$.  Finally, $\de < \si$, and $N_u^\si
\se N_u$ and so the conclusion is clear.

So, suppose we have defined  $(r_k:\ k \leq j)$
satisfying (a) - (c).  We define $r_{j + 1}$
and show that (a) - (c) are preserved.  Since
(b) clearly corresponds to (2),
and since we take $r = r_\si$, this will complete the proof,
once we show how (1) follows from (a).
This, however, is easy, since $(dom\ r) \sm (dom\ p^*) =
\bigcup \{(dom\ r_{k+1}) \sm (dom\ r_k):\ k < \si\}$,
and by (a), this last is indeed included in $\bcu \{N_u:\ u \in J\}$.

For (c) it will suffice to have $r_j \apr{\k_2} r_{j+1}$,
which will be clear from construction, as will the second part of
(b).  Thus, we must show that there is $q$ satisfying:

\roster
\item"{($\a$)}"\ \ $r_j \apr{\k_2} q$,
\item"{($\be $)}"\ \ $(dom\ q) \sm (dom\ r_j) \se N_{u_j}$,
\item"{($\ga $)}"\ \ if $card\ u_j = n$, then $q$
decides the value of $\bold c (u_j)$,
\item"{($\de $)}"\ \ for all
$u \in J,\ q\rest N_u \in N_u$.
\endroster

We first argue that it will suffice to find $q$ satisfying ($\a$) - ($\ga$),
since any such $q$ will automatically satisfy ($\de$).  For this,
note that if $q$ satisfies ($\a$), then $(dom\ q) \sm (dom\ r_j)$ has
power $< \th_{\k_2} = \k_1 \leq \h$.  Thus, for $u \in J,\
(q \sm r_j)\rest N_u$ is a subset of $N_u$ of power $< \h$ and therefore,
$(q \sm r_j) \in N_u$.  But $q\rest N_u = (q \sm r_j)\rest N_u
\cup r_j\rest N_u$,
and by induction hypothesis, (b), for $j,\ r_j\rest N_u \in N_u$.  The
conclusion is then clear.

To find $q$ satisfying ($\a$) - ($\ga$) is trivial if $card\ u_j < n$,
so assume $card\ u_j = n$.  Applying induction hypothesis (b), with
$k = j$ and $u = u_j$, we have $r_j\rest N_{u_j} \in N_{u_j}$.
Since the maximal antichain in $\bold Q^{apr}_{\k_2,\ p^*}$
deciding $\bold c (u_j)$ is a member of $N_{u_j}$,
and since $p^* \in N_{u_j}$, we easily
find $q^\prime \in N_{u_j}$ such that $p^* \cup r_j\rest N_{u_j} \apr{\k_2}
q^\prime$ and
such that $q^\prime$ decides the value of $\bold c (u_j)$.
Note that, again, since $(dom\ q^\prime) \sm (p^* \cup
(dom\ r_j\rest N_{u_j}))$ has
small cardinality, compared to the closure of $N_{u_j}$, we will
also have $(dom\ q^\prime) \sm (dom\ p^*) \se N_{u_j}$.
But this makes it clear that if
we take $q\ :=\ q^\prime \cup r_j$, then $q$ is as required.  This completes
the proof of Claim 2.
\enddemo

Now, let:

$$r^*\ :=\ p^* \cup \bigcup\{h_{u,\ v}((r \sm p^*)\rest N_u):\
u \in J,\ v \in \tilde J,\
u \sim v \}.$$

We will show that $r^* \in Q$
and that whenever $u \in J,\ v \in \tilde J$
and $u \sim v,\ p^* \cup h_{u,\ v}(r\rest N_u) \leq r^*$.
We first note that this
suffices for the proof of the Lemma in our special case,
since then clearly $r^*$ forces that $(B_i:\ i < \si)$
is as required.

The following is the heart of the matter, and is an easy
consequence of (7) of (2.1), and the arguments
for the first part of (2) of Claim 2, above.

\proclaim{Proposition}  Suppose that for $k = 1,\ 2,\
u_k \in J, v_k \in \tilde J$ and $u_k \sim v_k$.  Let
$N_k\ :=\ N_{u_k},\ N\ :=\ N_1 \cap N_2$ and let
$\tilde N = N_{u_1 \cap u_2}$ (so that, by (5) of (2.1),
$N \se \tilde N$).  Let $h_k\ :=\ h_{u_k,\ v_k}$.
Then, $(r \sm p^*)\rest N \in N\ \&\
h_1((r \sm p^*)\rest N) =  h_2((r \sm p^*)\rest N)$.
\endproclaim

\demo{Proof}  To see that $(r \sm p^*)\rest N \in N_k$, we argue as in the
proof of Claim 2:  $(r \sm p^*)\rest N$ is a subset of $N_k$ of
small cardinality compared to the closure of $N_k$.  But then,
since $h_1$ and $h_2$ are compatible functions, by (7) of
(2.1), the conclusion is clear.
\enddemo

\proclaim{Corollary}  $r^* \in Q$
and whenever $u \in J,\ v \in \tilde J$
and $u \sim v,\ p^* \cup h_{u,\ v}(r\rest N_u) \leq r^*$.
\endproclaim

\demo{Proof}  It is immediate from the Proposition, that
the $p^* \cup h_{u,\ v}((r \sm p^*))$ are pairwise
compatible as functions.  To complete the proof that
$r^* \in Q$, we must verify (a) of (1.4).  So, suppose
$i < \mu,\ \nu \in K$.

We consider separately the cases $\nu > \k_2,\ \nu < \k_2$,
and the hardest case, $\nu = \k_2$.  If $\nu > \k_2$, then
$\k_2 \leq \th_\nu$ and we taking the union of fewer than
$\th_\nu$ conditions, so there is no problem.  If $\nu < \k_2$,
then $\th_\nu < \k_1 \leq \h$, so for all $v \in \tilde J$,
either $\ci\nu \se N_v$ or $\ci\nu \cap N_v = \es$, and
then the conclusion is also easy.  So, suppose that $\nu = \k_2$,
i.e., $\th_\nu = \k_1$.
It is here that we use that $\k < \k_1$; this permits us to
argue as in the case where $\nu > \k_2$:  we are taking the
union of fewer than $\th_\nu$ conditions, and there is no
problem.

To complete the proof of the Corollary, we must see that
(c) of (1.4) holds (since (b) is clear).
So, once again, assume $\nu \in K$.  We must see that for all
$u \in J,\ v \in \tilde J$ such that
$u \sim v$, there are fewer than $\th_\nu$ many
$A \in \m /E_\nu$ such that $A$ grows from
$dom(p^* \cup h_{u,\ v}(r\rest N_u))$
to $dom\ r^*$.
Once again, is the proof that (1.4) (a), we consider separately
the cases $\nu  > \k_2,\ \nu < \k_2$ and $\nu = \k_2$.  Once
again, the hypothesis that $\k < \k_1$ allows us to assimilate
the case $\nu = \k_2$ to the case $\nu > \k_2$, since what is
really at issue is that we are taking the union of fewer
than $\th_\nu$ conditions, and as before, when $\nu = \k_2,\
\th_\nu = \k_1$.  In the remaining case, where $\nu < \k_2$,
once again we have that for all $i < \mu$ and all
$w \in \tilde J$, either $\ci\nu \se \N_w$ or $\ci\nu \cap N_w =
\es$, with the former holding if $\ci\nu \in N_w$.

So, suppose that $\nu < \k_2$ and fix such $u,\ v$, suppose $i < \mu$
and $\ci\nu$ grows from $dom(p^* \cup h_{u,\ v}(r\rest N_u))$
to $dom\ r^*$.  But then there are $t \in J,\ w \in \tilde J$
such that $t \sim w$ and $\ci\nu \cap dom(p^* \cup h_{t,\ w}(r\rest N_t))
\not\se dom(p^* \cup h_{u,\ v}(r\rest N_u))$.  Then, $\ci\nu \in
N_v \cap N_w$.  But then $\ci \se N_v \cap N_w$.  Therefore,
letting $b \in N_u,\ c \in N_t$ be such that $\ci\nu = h_{u,\ v}(b) =
h_{t,\ w}(c)$, we clearly have $b = c$ and $b \se N_u \cap N_t$.
But then, letting $x\ :=\
(dom\ r \sm dom\ p^*) \cap b,\ x \in N_u \cap N_t$, since,
once again, $x$ is a subset of each, small in cardinality
compared to the closure of each.  So $h_{u,\ v}(x) = h_{t,\ w}(x)$,
but this is a contradiction, since then,
$\ci\nu \cap dom(p^* \cup h_{u,\ v}(r\rest N_u)) =
((\ci\nu \cap dom\ p^*) \cup h_{u,\ v}(x)) =
((\ci\nu \cap dom\ p^*) \cup h_{t,\ w}(x)) =
\ci\nu \cap dom(p^* \cup h_{t,\ w}(r\rest N_t))$.  This completes the
proof of the Corollary, and therefore of the Lemma, in our special case.
\enddemo
\1\1
To handle the case $\tau = \k_2$, we take $B\ :=\
\bigcup\{B_i:\ i < \si\}$, we replace $D$, above, by
$D^\prm\ :=\ \{ a \in [B]^n:\ card(a \cap B_i) \leq 1\text{, for all }
i < \si\}$, and we take our
$\pr {\k_2}\text{-increasing}$ sequence from $Q,
\ \orh{p} = (p_j: j < \eta)$,
to satisfy:
\1\1

\roster
\item"{(1*)}"  $\eta \leq \k$, and $p_0 = p$,

\item"{(2*)}"  for each $\orh{\a} = (\a_1,\ \cdots,\ \a_n) \in D^\prm$, there is
$j = j_{\orh{\a}}$ such that in $\Qapr{\k_2,\ p_{j+1}}$, there is a predense
set, $I_{\orh{\a}}$ of conditions deciding $\bold c (\orh{\a})$.
\endroster
\1\1

Now $\eta \leq \k \leq \k_1 < \tau$, and then the rest of the proof goes
through easily, as above.




\subheading {(2.3)  THE REMAINING CASE}

We show how to handle the case $\si = \la$, under
the additional hypothesis that $\la$ is not strongly
inaccessible.
Our approach to this case involves building a condition
$r$, as in Claim 2 of (2.2).
Recall our choice of $(\ga_j:\ j < \si )$ and our definitions of
$J$ and $\tilde J$ from (2.2).
For each $u \in J$, with $card\ u = n$, we shall
build a condition $r_u$.  Letting $u^*\ :=\ \{\ga_1,\ \cdots ,\ \ga_n\}$,
for $u \in J$ with $card\ u = n$, we shall take
$r_u$ to be the isomorphic copy, by $h_{u^*,\ u}$, of a
condition $p_u$.
We will want to essentially have $p_u \sse N_{u^*}$, so that
this makes sense.  We shall associate to $u$
an appropriate $n\text{-tuple}$ of branches through a
certain tree, and we will build $p_u$ as a limit through
these branches of conditions indexed by the nodes of the tree.
The main work will be to construct these conditions indexed
by the nodes of the tree.

The properties of the $h_{u,\ v}$ can then be invoked,
just as in (2.2), to
guarantee that the union of the $r_u$ is in fact a condition
extending each $r_u$.  Since there is nothing essentially
new here, compared to the situation in (2.2), we shall
not repeat the argument.  The union of the $r_u$ will then
be our $r$, and once we have $r$,
the proof proceeds exactly as in
(2.2).  We first describe the tree.
\1

\subheading {(2.3.1)  The tree}

Recall that we are assuming that for some $\de < \la,\ 2^\de \geq \la$.
This guarantees the existence of a cardinal $\de < \la$ and a tree,
$T \sse \bigcup\{\ ^\ga 2:\ \ga < \de\}$ and satisfying
\1
\roster
\item $card\ T < \la$,
\item $T$ has at least $\la$ branches (of length $\de$),
\item for all $\a < \de$, there is unique $\eta_\a \in T_\a$ such
that $\eta_\a \cup \{(\a,\ 1)\} \in T$.
\endroster
\1
In what follows,
conventionally, we shall take $T_\de$
to be the set of all branches of length $\de$
through $T$.  Also,
if $\a \leq \de$, we shall use
$[T_\a]^n$ denote the set of $n\text{-element}$
subsets of $T_\a$, though this might seem to
conflict with the notation $[S]$ for the set
of all maximal branches through the tree $S$.
For $\a \leq \de$, we view the elements of $[T_\a]^n$ as
sequences, $\orh\eta = (\eta_1,\ \cdots,\ \eta_n)$ which
are increasing for lexicographic order.
For such $\a$ and $\orh\eta$ and $\be < \a$, we let
$\orh\eta |\be\ :=\ (\eta_1|\be,\ \cdots,\ \eta_n|\be )$.

\1
\subheading {(2.3.2)  The construction of $r$}

The argument has many features in common with that the proof
of Claim 2, in (2.2).
We work by recursion
on $\a \leq \de$ to define conditions, $p_{\orh\eta} \in N_{u^*}$,
for $\orh\eta \in [T_\a]^n$,
(so for $\a = \de$, the $\orh\eta$ are $n-\text{tuples}$ of $\de-
\text{branches}$ through $T$) which satisfy the following properties:
\1
\roster
\item"{(a)}"  if $\a < \be \leq \de,\ \orh\eta \in [T_\a]^n,\
\orh\eta^\prm \in [T_\be]^n$ and for $1 \leq i \leq n,\
\eta_i \sse \eta^\prm_i$, then $p_{\orh\eta} \leq p_{\orh\eta^\prm}$,

\item"{(b)}"  if $\a \leq \de$ is a limit ordinal and $\orh\eta \in [T_\a]^n$,
then $p_{\orh\eta}$ is the lub in $\bold Q$ of the $p_{\orh\eta |\be}$,
for $\be < \a$ such that the $\orh\eta_i|\be$ are all distinct, for
$1 \leq i \leq n$,

\item"{(c)}"  if $\orh\eta,\ \orh\eta^\prm \in [T_\a]^n$, we let
$w\ :=\ \{ i:\ \eta_i = \eta^\prm_i\}$ and we require that
$p_{\orh\eta}|N_{\{ \ga_i:\ i \in w\} } =
p_{\orh\eta^\prm}|N_{\{ \ga_i:\ i \in w\} }$,

\item"{(d)}"  each $p_{\orh\eta}$ decides the value of $\bold c (u^*)$.
\endroster

At limit stages $\a \leq \de$ and
$\orh\eta^* \in [T_\a]^n$, let $\be^*$ be the least $\be < \a$ such that
the projection of $\orh\eta^*$ to level $\be$ is in $[T_\be]^n$, i.e.,
such that the restrictions to $\be$ of the terms of $\orh\eta^*$
are pairwise distinct, and for $\be < \ga < \a$, let
$\orh\eta^*|\ga$ be the projection of $\orh\eta^*$ to level $\ga$.
We then define $p_{\orh\eta^*}$ to be the union of the
$p_{\orh\eta^*|\ga}$ for $\be \leq < \a$.  Since $\a \leq \de < \la$,
this presents no problem.

So, consider a successor stage, $\a + 1$.
Note that in this case, we are {\bf below } $\de$, so our $\orh\eta$
are $n-\text{tuples}$ of actual elements of $T$.
For $i = 0,\ 1$, we let $\eta_\a^i$ denote $\eta_\a \cup \{ (\a,\ i)\}$,
so $\eta_\a^0,\ \eta_\a^1 \in T_{\a + 1}$

The main complication arises in dealing with $\orh\eta \in [T_{\a + 1}]^n$
such that $\eta_\a^0,\ \eta_\a^1$ are both terms of $\orh\eta$, since then
there is no preexisting \lq\lq projection" of $\orh\eta$ to level $\a$.
Such $\orh\eta$ are called {\it new}.
We let $(\orh\eta^j:\ j < j^*)$ enumerate the new $\orh\eta$, where
$j^*$ is a cardinal, which, in view of (1) of (2.3.1),
and the fact that we are below $\de$, must be less than $\la$.

We first define, essentially
trivially, for {\bf all } $\orh\eta \in [T_{\a + 1}]^n$,
conditions $p^0_{\orh\eta}$ which satisfy (a) and (c) above.  If $\orh\eta$
is new, then there is no reason to believe that $p^0_{\orh\eta}$ satisfies
(d), however.  Then, we define, for {\bf all } $\orh\eta$,
by recursion on $j \leq j^*$, an increasing sequence of conditions
$(p^j_{\orh\eta}:\ j \leq j^*)$ (with
each $p^j_{\orh\eta} \in N_{u^*}$.  At limit stages, $j$, for each $\orh\eta$,
we take $p^j_{\orh\eta}$ to be the union of the $p^i_{\orh\eta}$, for
$i < j$; since $j^* < \la$, this creates no difficulties.  To define
$p^{j + 1}_{\orh\eta}$, we first extend $p^j_{\orh\eta^j}$ to a condition
$q^j_{\orh\eta^j} \in N_{u^*}$ which decides $\bold c (u^*)$; for
$\orh\eta \neq \orh\eta^j$, we take $q^j_{\orh\eta} = p^j_{\orh\eta}$.
Finally, we extend each of the $q^j_{\orh\eta}$ to $p^{j + 1}_{\orh\eta}$
in such a way as to recover (c), above.
This completes the construction of the $p_{\orh\eta}$, and in
particular of the $p_{\orh\eta}$ for the $\orh\eta \in [T_\de]^n$.

Enumerate as $(\eta_i:\ i < \la)$, without repetitions, some $X \sse
T_\de$ with $card\ X = \la$.  Finally, for $u =
\{ \ga_{i_1},\ \cdots,\ \ga_{i_n}\} \in J$, with $i_1 < \cdots < i_n$, let
$\orh\eta = (\eta_{i_1},\ \cdots,\ \eta_{i_n})$ and define $p_u$ to be
$p_{\orh\eta}$.  This completes the proof.

\enddocument
\end